\documentclass[12pt]{article}
\usepackage{amsmath}
\usepackage{amssymb}

\newcommand{\re}{\mathbb R}
\newcommand{\z}{\mathbb Z}
\newcommand{\bs}{\cal B}
\newcommand{\msp}{(S,\bs,{\mu})}
\newcommand{\bp}{(S,\bs)}
\newcommand{\lps}{L^{\alpha}(S,\bs,\mu)}
\newcommand{\ma}{{\cal S}{\alpha}{\cal S}}
\newcommand{\str}{\{X_t,t \in \re\}}
\newcommand{\stp}{\{X_t,t \in T\}}
\newcommand{\fsp}{\{f_t,t \in T\}}
\newcommand{\sps}{\{X_t,t \in T_0 \}}
\newcommand{\hn}{H\neq 1/\alpha}
\newcommand{\ha}{H-1/\alpha}
\newcommand{\sug}{Surgailis et al (1993)}
\newcommand{\fn}{\longmapsto}
\newcommand{\sm}{\stackrel{d}{=}}

\begin{document}
\begin{titlepage}

\title{On the mixing structure of stationary increment and self-similar $\ma$ processes}

\author{Donatas Surgailis\footnote{Institute of Mathematics and Informatics, Lithuanian Academi of Science, 2600 Vilnius, Lithuania} \and Jan Rosinski\footnote{Department of Mathematics, University of Tennessee, Knoxville, TN-37996} \and V. Mandrekar\footnote{Department of Statistics and Probability, Michigan State University, E. Lansing, MI-48824} \and \fbox{Stamatis Cambanis}}
\date{1998}

\maketitle

\begin{abstract}
 Mixed moving average processes appear in the ergodic decomposition of stationary symmetric $\alpha$-stable  ($\ma$ )  processes. They correspond to the dissipative part of ``deterministic'' flows generating  $\ma$ processes (Rosinski, 1995).  Along these lines we study stationary increment and self-similar $\ma$ processes.  Since the classes of stationary increment and self-similar  processes can be embedded into the class of stationary processes by the Masani and Lamperti transformations, respectively, we characterize these  classes of  $\ma$ processes in terms of nonsingular flows and the related cocycles.  We  illustrate this approach considering   various  examples of self-similar mixed moving average $\ma$ processes introduced in (Surgailis, Rosinski, Mandrekar and Cambanis, 1992). 
 \bigskip \bigskip
 
{\footnotesize 

{\sl AMS} 1991 {\sl subject classifications. } Primary 60G10; secondary 60G07, 60E07, 60G57. 
 
}
\end{abstract}

\end{titlepage}
\section*{0\hspace{.2in} Introduction}
The mixed moving average $\ma$ processes were introduced by the authors (1993). Subsequently, J. Rosinski (1995) showed that these form or the dissipative part in the ergodic decomposition of stationary $\ma$ process. In the work of the authors (1992), mixed moving average, self-similar stationary increment processes were considered and shown to provide a large class of self-similar processes with stationary increment well beyond the linear fractional stable motions.

In the first part of this work we consider the ergodic decomposition of stationary increment and self-similar processes separately. In both cases, we associate a stationary $\ma$ process and use Rosinski's decomposition. Although in self-similar case, one can use Lamperti's construction to associate a stationary $\ma$ process, one needs to do some technical work to extend Masani's work (1976) for stationary increment case for all $\alpha >0 \;(\alpha <1,$ in particular). This is shown in the Appendix and Theorem 2.2. First, we look at dissipative and conservative parts for $\ma$ processes, which are both self-similar and with stationary increments. The processes studied by the authors (1992) which include classical examples are shown to be a subclass of dissipative processes studied here and finally, we study a large class of conservative $\ma$ stationary increment self-similar processes. These turn out to be rotationally mixing.

\section{Preliminaries and Notation}
A stochastic process $\stp$ is said to be symmetric $\alpha$ stable $(\ma)$ if any finite linear combination $\sum a_i X_{t_i},\; a_i \in \re$ and $t_i \in T$ has $\ma$  distribution. A family of functions $\fsp \subseteq \lps$ where $\msp$ is a standard Borel space with $\mu \;\; \sigma$-finite ([1]) is said to be a spectral representation of a $\ma$ processes $\stp$ if
$$\stp \sm \left\{ \int_{\cal S}f_t(u)dM(u), t \in T\right\} \leqno(1.1)$$
where $M$ is a independently scattered random measure on $\bs$ such that
$$E \exp\{it M(A)\}= \exp\{ -|t|^{\alpha}\mu(A)\},\qquad t \in \re,\, A \in \bs .$$
We also consider complex stable  processes. However in the complex case we restrict our attention to those $\stp$ so that $\sum a_i X_{t_i},\; a_i \in C,\;t_i \in T$ are rotationally invariant stable distributions. In this case a family of complex $\alpha$ -integrable function $\fsp$ defined on Borel space $\msp$ ($\mu \;\; \sigma$-finite) is called a spectral representation of $\stp$ if (1.1) holds with a complex independently scattered measure $M$ satisfying
$$E \exp\{i\,Re (\overline{tM(A)})\}= \exp\{ -|t|^{\alpha}\mu(A)\}, \qquad t \in \re, A \in \bs.$$
A stochastic process $\stp$ is called separable in a probability if there exists a countable set $T_0 \subseteq T$ such that the set of random variables $\sps$ is dense in the set $\stp$ with topology of convergence in probability. Every separable in probability $\msp$ process has a spectral representation with $S$ unit interval and $\mu$ Lebesgue measure on $S$ (see Kuelbs (1973) and Hardin (1982) for complex case). Conversely, if $\msp$ has a spectral representation defined on standard Borel space with a $\mu \; \sigma$-finite, then it is separable in probability. A spectral representation $\fsp$ (for $T$ separable metric space) is measurable if $(s,t) \fn f_t(s)$ is measurable with respect to $\sigma$-algebra of $S\times T$. Every measurable $\ma$ process has measurable spectral representation (Rosinski and  Woyczynski (1986)). Also every measurable process is  separable in probability. Throughout we consider only measurable processes. We also identify equivalent represen!
 tation in the sense of Rosinski (1994).

Let $\{f_t\}_{t \in T} \subseteq \lps$ be a collection of functions and $\sigma_r(\fsp)$ denote the smallest $\sigma$-field generated by extended-valued functions $f_t/f_{\tau} \;\; t,\tau \in T$. Following Hardin (1982) we give the following definition.

{\bf 1.2 Definition.} A spectral representation $\fsp \subseteq \lps$ of $\msp$ process is said to be minimal if $supp \fsp = S \;\; \mu- a.e.$ and for every $B\in \bs$ there exists an $A \in \sigma_r(\fsp)$ such that $\mu(B\Delta A)=0$.

The following theorem is due to Hardin(1982).

{\bf 1.3 Theorem. } {\em Every separable in probability $\ma$ process has a minimal spectral representation. One can choose $S$ as a unit interval or a countable discrete set or the union of the two and $\mu$ as the direct sum of Lebesgue measure on unit interval and counting measure on the discrete set.}

In this work we shall study representation for certain class of $\ma$ processes by relating them to the representation of stationary $\ma$ processes studied by Rosinski (1995). For this, we need the following concepts. Let $T=\re$ on $\z$. A stochastic process $\stp$ is said to be stationary if for every $\tau \in T$,
$$\{X_{t+\tau},t \in T\} \sm \stp.$$
A family $\{\varphi_t,t \in T\}$ of measurable maps on a Borel space $\bp$ is said to be a flow on $S$ if, for $t_1,t_2 \in T$ and $s\in S$,
$$\varphi_{t_1+t_2}(s)=\varphi_{t_1}(\varphi_{t_2}(s)) \;\;\;\mbox{and} \;\;\; \varphi_0(s)=s.$$
A flow $\{\varphi_t,t \in T\}$ is measurable if $(t,s) \fn \varphi_t(s)$ is measurable. Given  a $\sigma$-finite measure $\mu$ on $\bp$, $\{\varphi_t,t \in T\}$ is said to be $\mu$ non-singular if
$$\mu(\varphi_t^{-1}(A))=0 \;\; \mbox{iff} \;\; \mu(A)=0 \;\; \mbox{for} \;\; t\in T, A \in \bs .$$
We denote by $\rho_t(s)=\frac{d\mu_0\varphi_t}{d\mu} (s)$. Then
$$\rho_{t+u}(s)=\rho_t(s)\rho_u(\rho_t(s))=\rho_u(s)\rho_t(\rho_u(s)).\leqno(1.4).$$
A measurable map $a_t(s)$ from $(T \times S) \fn G$, a second countable group is called a cocycle for a measurable flow $\{\varphi_t,t \in T\}$ if for $t_1,t_2 \in T$
$$a_{t_1+t_2}(s)=a_{t_2}(s)a_{t_1}(\varphi_{t_2}(s)) \qquad a.e. \mu .\leqno(1.5)$$

A cocycle $\{a_t,t \in T\}$ is said to be a coboundary if there exists a measurable $b:S \fn G$ so that $a_t(s)=b(a_t(s))b^{-1}(s) \;\; a.e. \;\mu$ for each $t$. The following Theorem is due to Rosinski(1995).

{\bf 1.6 Theorem } {\em Let $\fsp \subseteq \lps$ be a measurable minimal spectral representation of a measurable stationary $\ma$ process $\stp$. Then there exist a unique (modulo $\mu$) non-singular flow $\{\varphi_t,t \in T\}$ on $({\cal S},\mu)$ and a cocycle $\{a_t,t \in T\}$ for $\{\varphi_t,t \in T\}$ taking values in $\{-1,1\}\;(\{|z|=1\}$ in complex case) so that for each $t \in T$
$$f_t=a_t \rho_t^{1/\alpha}f_0\circ\varphi_t \qquad \mu \;a.e. \leqno(1.7)$$}

\section{Spectral Representation of Processes With Stationary increments}
A stochastic process $\stp \;(T=\re,\z)$ is said to have stationary increments (s.i.) if for all $h\in T$
$$\{X_{t+h}-X_{s+h},\;t,s \in T\}\sm \{X_t-X_s,\;t,s \in T\}$$
We observe that

{\bf 2.1 Proposition: } {\em Every measurable s.i. process $\stp$ is continuous.}

{\it Proof: } Consider the F-norm on $L^0(\Omega,P)$ given by $\|y\|=E\min (|y|,1)$ and define for $\epsilon >0$ (fixed),
$$B_t=\{s\in \re:\|X_t-X_s\|<\epsilon\}.$$
Under the assumption and a result by Cohn(1972), the map
$$t\in \re \fn X_t\in L^0(\Omega,P)$$
is Borel and has separable range. Thus we can choose a sequence $\{t_n\}\subseteq \re$ so that $\{B_{t_n}\}$ are Borel and $\re=\cup_nB_{t_n}$. Hence, there exits at least one $t_n$ so that Lebesgue measure of $B_{t_n}$ is positive. By Steinhaus Lemma, the set $B=B_{t_n}-B_{t_n}$ contains an open interval centered at zero, say, $(-\delta,\delta)$. If $|s-t|<\delta$, then $s-t=u-v\;\;(u,v, \in B_{t_n})$ and
$$\|X_t-X_s\|=\|X_{s-t}-X_0\|=\|X_u-X_v\|\leq \|X_u-X_{t_n}\|+ \|X_v-X_{t_n}\| <2 \epsilon.$$
This proves the uniform stochastic continuity of $X.\;\;\square$

We now give an extension of a main result of Masani(1976) in the form used by Cambanis and Maejima (1989). Since every stochastic processes can be viewed as curve in $L^0(\Omega,P)$, the condition $X\in {\cal R}_m([a,b]),L^0(\Omega,P))$ makes sense (see Appendix). This amounts to verifying whether or not the integral
$$\int_a^b \varphi(t)X_tdt$$
exists as the limit in probability of Riemann sums, for every $\varphi \in D[a,b]$. Clearly, if sample paths of $X$ are Riemann integrable on $[a,b]$, then $X\in {\cal R}_m([a,b]),L^0(\Omega,P))$. If $X$ is s.i. process, then the so is $X'_t=X_t-X_0$ and vice-versa. Thus, without loss of generality, we may assume $X_0=0$.

{\bf 2.2 Theorem: } {\em Let $X$ be an s.i. measurable process with $X_0=0$ such that
$$
\begin{array}{ll}
(i) & X\in {\cal R}_m([0,1]),L^0(\Omega,P))\\
(ii) & E\log^+|\int_0^1 e^t X_tdt|<\infty,\\
(iii) & E\log^+|X_1|<\infty.
\end{array}$$
Then
$$Y_t=X_t-\int_{-\infty}^t e^{-(t-s)}X_sds \leqno(2.3)$$
is a well-defined stationary stochastically continuous process so that, for every $s<t,\; Y\in {\cal R}_m([s,t]),L^0(\Omega,P))$ and
$$X_t=Y_t-Y_0+ \int_0^tY_udu \qquad t\in \re. \leqno(2.4)$$
Conversely, if $Y\in {\cal R}_m([s,t]),L^0(\Omega,P))$, stationary then (2.4) defines and s.i. processes $X_t$ with $X_0=0$.}

{\it Proof: } First we show that
$$I=\int_{-\infty}^0 e^uX_udu\;:=\lim_{n\to -\infty}\int_{-n}^0 e^uX_udu$$
exists with probability one. Consider
$$Z_n=\int_{-n}^{-n+1}e^uX_udu, \qquad n \geq 1.$$
Since $\{X_{v-n}-X_{-n}\;:v \in [0,1]\} \sm \{X_v\;: v \in [0,1]\}, \;\;Z_n$ is well defined and
$$Z_n=\int_0^1e^{v-n}X_{v-n}dv=e^{-n}V_n+(e-1)e^{-n}X_{-n},$$
where $V_n=\int_0^1e^v[X_{v-n}-X_{-n}]dv \sm \int_0^1e^vX_vdv=V_0$.

First we will show that $\sum e^{-n}|V_n|<\infty$ a.s. This follows from (ii) since
\begin{eqnarray*}
\sum P\{e^{-n}|V_n|>e^{-n/2}\} & = & \sum P\{|V_0|>e^{n/2}\}\\
& = & \sum P\{2 \log^+|V_0|>n\}<\infty
\end{eqnarray*}
Next we show that $\sum e^{-n}|X_n|< \infty$ a.s. Put $W_k=X_{-k}-X_{-k+1};\;X_{-n}=\sum_{k=1}^nW_k$ and $W_k\sm-X_1$. Hence
\begin{eqnarray*}
\sum_{k=1}^{\infty}e^{-n}|X_{-n}| & \leq & \sum_{k=1}^{\infty}\sum_{n=k}^{\infty} e^{-n}|W_k|\\
& = & (1-e^{-1})^{-1}\sum_{n=k}^{\infty} e^{-k}|W_k| < \infty
\end{eqnarray*}
by (iii) and the same argument as above. Thus $\sum|Z_n| < \infty$ a.s. and $I$ is well defined. This fact enables us to formally define
$$\int_{-\infty}^t e^{t-s}X_sds\;:=e^{-t}I+e^{-t}\int_{-\infty}^t X_sds.$$
Thus $Y$ is well-defined in (i) and, by Lemma 1 of Appendix, $Y\in {\cal R}_m([s,t]),L^0(\Omega,P))$ for every $s<t$. Continuity of Riemann's integral implies the stochastic continuity of $Y$ and since
$$Y_t=\int_{-\infty}^t e^{t-s}[X_t-X_s]ds$$
and $X$ has s.i., $Y$ is a stationary process. Now, using Lemma 1 of the Appendix, we get
\begin{eqnarray*}
\int_0^t Y_sds & = & \int_0^t X_sds - \int_0^t\left [e^{-s}I+e^{-s}\int_0^s e^u X_udu\right ]ds\\
& = & \int_0^t X_sds +(e^{-t}-1)I-\int_0^t\left (\int_u^se^{-s}ds\right ) e^u X_udu\\
& = & e^{-t}I+\int_0^t e^{t-s}X_udu-I\\
& = & X_t-Y_t+Y_0.
\end{eqnarray*}
This proves (2.4). Since the converse is obvious, the proof of Theorem 2.2 is complete. $\square$

Let us now assume that $\stp$ is a measurable $\ma$ process so that $X\in {\cal R}_m([a,b]),L^0(\Omega,P))$. In view of proposition 2.1, for $\alpha >1$, we get that this assumption is trivially valid. Using Theorem 2.2, we get
$$X_t-X_0=Y_t-Y_0-\int_0^t Y_udu \leqno(2.5)$$
where $\{Y_t, \; t \in \re\}$ is stochastically continuous stationary process given in (2.4). We note that $\{Y_t, \; t \in \re\}$ is a measurable $\ma$ process.

We shall denote by $\|\cdot\|_{\alpha}$ the Schilder norm.

Let us assume that $\stp$ has a spectral representation
$$X_t= \int_S f_t(s)dM(s) \leqno(2.6)$$
In view of (2.5) condition for minimality  for $\{Y_t, \; t \in \re\}$ can be expressed in terms of $\fsp$.

{\bf 2.7 Theorem: }{\em A measurable $\ma$ process $\stp$ has s.i. iff there exists a spectral representation (2.6) and a group $\{P^t, t \in \re\}$ of isometrics of $\lps$ so that
$$P^{\tau}(f_t(\cdot)-f_0(\cdot))=f_{t+\tau}(\cdot)-f_{\tau}(\cdot).$$}

The proof is simple, using
$$\|(f_t(\cdot)-f_0(\cdot\|_{L^{\alpha}(\mu)} = \|X_t-X_0\|_{\alpha}.$$
Under additional conditions, we can give precise form of $f_t$ and describe the isometries $P^{\tau}$ using Theorem 1.6.

{\bf 2.8 Theorem :} {\em Let $\stp$ be a measurable $\ma$ s.i. process such that\\$X\in {\cal R}_m([a,b]),L^0(\Omega,P))$ with associated stationary process having minimal representation. Then there exists a unique (modulo $\mu$) non-singular flow $\{\varphi_t,t \in T\}$ on $\bp$, a cocycle $\{a_t,t \in T\}$ for $\{\varphi_t,t \in T\}$ taking values in $\{-1,1\}\;(|z|=1$ in complex case) and a function $g_0\in \lps$ such that for $t\in T$
$$f_t(s)-f_0(s)=a_t\rho_t^{1/\alpha}g_0\circ\varphi_t(s)-g_0(s) - \int_0^ta_u\rho_u^{1/\alpha}g_0\circ\varphi_u(s)du \leqno a)$$
and for $t,t' \in T$ and each $\tau \in T$
$$P^{\tau}(f_t(\cdot)-f_{t'}(\cdot))= a_{\tau}\rho_{\tau}^{1/\alpha}[(f_t \circ\varphi_{\tau})-(f_{t'} \circ\varphi_{\tau})]\leqno b)$$}

{\it Proof: } From (2.5) and Theorem 1.6 we get
\begin{eqnarray*}
\lefteqn{\int_S(f_t-f_0)dM(s) =} \\
& & \qquad \int_S a_t\rho_t^{1/\alpha}g_0\circ\varphi_tdM(s)-\int_S g_0dM(s)-\int_0^t\int_S a_u\rho_u^{1/\alpha}g_0\circ\varphi_udM(s)
\end{eqnarray*}
where $Y_t=\int  a_t\rho_t^{1/\alpha}g_0\circ\varphi_tdM(s)$. Using Samorodnitsky (1992) and Rosinski (1994) we get (a). To obtain part (b), we observe $t'<t \in T$ and $ \tau\in T$
\begin{eqnarray*}
\lefteqn{P^{\tau}(f_t-f_{t'})}\\
&&\qquad =f_{t+\tau}-f_{t'+\tau}\\
&&\qquad =a_{t+\tau}\rho_{t+\tau}^{1/\alpha}g_0\circ\varphi_{t+\tau}-\int_{t'+\tau}^{t+\tau} a_u\rho_u^{1/\alpha}g_0\circ\varphi_udu
-a_{t'+\tau}\rho_{t'+\tau}^{1/\alpha}g_0 \circ\varphi_{t'+\tau}
\end{eqnarray*}
Note
$$\int_{t'+\tau}^{t+\tau} a_u\rho_u^{1/\alpha}g_0\circ\varphi_udu = \int_{t'}^t a_{u+\tau}\rho_{u+\tau}^{1/\alpha}g_0\circ\varphi_{u+\tau}du.$$
Now use $\varphi_{t_1+t_2}=\varphi_{t_1}(\varphi_{t_2})$, (1.4) and (1.5) to get
$$f_{t+\tau}-f_{\tau}=a_{\tau}\rho_{\tau}^{1/\alpha}[(f_t \circ\varphi_{\tau})-(f_{t'} \circ\varphi_{\tau})].$$

{\bf 2.9 Corollary: } {\em If the measurable $\ma$ s.i. process $\stp$ is of the form $X_t=Y_t+Z$ where $\{Y_t,Z\}$ is an $\ma$ process with $\{Y_t,t\in T\}$ stationary, then $Z=-X_0-Y_0$ and
$$f_t(s)-f_0(s)=a_t\rho_t^{1/\alpha}g_0\circ\varphi_t(s)-g_0(s)$$}

{\bf 2.10 Corollary: } {\em If the measurable $\ma$ s.i. process is of the form $X_t=\int_0^t Y_udu$ with $\{Y_t,t\in T\}$ stationary measurable $\ma$ process then
$$f_t(s)=\int_0^t a_u\rho_u^{1/\alpha}g_0\circ\varphi_u du.$$}

Let us define, with $\lambda$ Haar measure on $T$
$$C=\{s\in S:\;\int_T|g_0(\varphi_t(s))|^{\alpha}\rho_t(s)\lambda(dt)<\infty \}$$
and
$$D=\{s\in S:\;\int_T|g_0(\varphi_t(s))|^{\alpha}\rho_t(s)\lambda(dt)=\infty \}.$$
Then $C\cup D=S$ is the Hopf Decomposition of the non-singular flow $\{\varphi_t,t\in T\}$ by Theorem 4.1 of Rosinski (1995).

We say that a measurable $\ma$ s.i. process is generated by a non-singular measurable flow $\{\varphi_t,t\in T\}$ on $\bp$ if
$$f_t(s)-f_0(s)=a_t\rho_t^{1/\alpha}g_0\circ\varphi_t(s)-g_0(s) - \int_0^ta_u\rho_u^{1/\alpha}g_0\circ\varphi_u(s)du $$
where $f_t \in L^{\alpha}\bp$ and
$$S=supp\{f_t(s)-f_0(s),\;t\in T\}=supp\{f_t(s)-f_{t'}(s),\;t,t'\in T\}$$
In view of Masani (1976) Theorem 3.15 and Theorem 2.2 we get that
$$S=supp\{g_0\circ\varphi_t(s),\;t\in T\}.$$
We have shown that any measurable $\ma$ s.i. process is generated by a measurable non-singular flow on a standard Borel Space $S$. Define
$$X_t^D=\int_Df_tdM,$$
$$X_t^C=\int_Cf_tdM.$$
Then $X_t-X_0\sm (X_t^D-X_0^D)+(X_t^C-X_0^C)$. Since the associated stationary processes are independent we get $(X_t^D-X_0^D)$ is independent of $(X_t^C-X_0^C)$. We also get upto a constant random variable the decomposition is unique in distribution.

Using Theorem 2.8 and Theorem 4.4 of Rosinski (1995) we get

{\bf 2.11 Theorem: } {\em Let $\stp$ be a measurable $\ma$ process generated by a dissipative flow. Then there exists a Borel space $W$, a $\sigma$-finite measure $Q$ on $W$, and a function $g\in L^{\alpha}(W\times T,Q \otimes\lambda)$ and an independently scattered $\ma$ measure $N$ with control measure $Q\otimes \lambda$, so that
$$\{X_t-X_0\}\sm\left\{\int_W \int_T\left [g(x,t+u)-g(x,u)-\int_0^t g(x,s+u)ds\right ]N(dx,du)\right\}.$$}
In particular by corollary 2.6 with $X_0=0$, we get
$$\stp \sm \left\{\int_W \int_T\left [g(x,t+u)-g(x,u)\right ]N(dx,du),\;t\in T\right\}$$

{\bf Remark: } Note that the distribution of $N(A,-B)$ is the same as $N(A,B)$ giving the fact $\stp \sm \{\int_W \int_T [g(x,t+u)-g(x,u)]N(dx,du),\;t\in T\}$. These processes were originally considered in Surgailis, et all (1992).

\section{Spectral Representation of self-similar Processes.}
Recall that a process $\str$ is self-similar with parameter $H>0\;(H$-ss) if $\{X(\lambda t), t\in \re\}\sm \{\lambda^HX(t),t\in \re\}$ for all $\lambda>0$. Let us assume that $\stp$ is $\ma$ process with representation
$$X(t)=\int_Sf_t(s)dM(s)\qquad \{f_t,t\in \re\}\subseteq L^{\alpha}\bp.$$

{\bf 3.1 Theorem: } {\em A separable in probability $\ma$ process $\str$ is $H$-ss iff there exists a spectral representation $\{f_t,t\in \re\}$ and a group $\{R^{\lambda},\lambda >0\}$ of isometries of $L^{\alpha}({\cal S},\mu)$ (i.e. $R^{\lambda_1 + \lambda_2}=R^{\lambda_1}R^{\lambda_2}$) such that
$$f(\lambda t,\cdot)=\lambda^H R^{\lambda}f(t,\cdot)\qquad\mbox{for } \lambda>0\mbox{ and } t\in\re.$$}

The proof follows by the definition.

Let us now recall that the Lamperti transformation which says that $\{X_t,t>0\}$ is $\ma\;H$-ss process if $\{Y(t)=e^{-Ht}X(e^t),t\in\re\}$ is stationary $\ma$ process. In particular, if $\{Y(t),t\in\re\}$
 stationary $\ma$ then $X(t)=t^{-H}Y(\log t)$ is $\ma \;H$-ss process.

{\bf 3.2 Theorem: } {\em Let $\{X_t,t\in\re_+\}$ be a measurable $\ma \;H$-ss process, separable in probability such that  $\{Y(t),t\in\re\}$ has minimal representation. Then there exists a unique (modulo $\mu$) non-singular flow $\{\tilde{\varphi}_t,t\in\re\}$ on $\bp$, a cocycle $\{a_t,t\in\re\}$ for $\{\tilde{\varphi}_t,t\in\re\}$ taking values in $\{-1,1\}\;(|z|=1$ in complex case) and a function $g_0\in L^{\alpha}\bp$ so that for $t>0$
$$f_t(s)=t^Ha_{\log t}\rho_{\log t}^{1/\alpha}g_0\circ \tilde{\varphi}_{\log t}(s)\leqno(3.3)$$
$$R^{\lambda}f_t(s)=a_{\log \lambda}(\rho_{\log \lambda})^{1/\alpha}f_t\circ \tilde{\varphi}_{\log \lambda},\qquad\lambda>0.\leqno(3.4)$$}

{\it Proof: } Using Lamperti transformation, we get for $t\in\re,\;\{Y(t)=e^{-Ht}X(e^t)=\int_Sg_t(s)dM(s)\}$ is stationary $\ma$ satisfying assumptions of Theorem 3.1 of Rosinski (1995) giving $g_t(s)=a_t\rho_t^{1/\alpha}g_0\circ \tilde{\varphi}_t(s)$ for $t \in \re$. We have $g_t(s)=e^{-Ht}f_{e^t}(s)$ or equivalently $t>0$, $f_t(s)=t^Hg_{\log t}(s)$. Using the above representation we get a). To obtain b) we observe that
$$f_{\lambda t}(s)=(\lambda t)^H a_{\log \lambda +\log t}(\rho_{\log \lambda +\log t})^{1/\alpha}f_1\circ\tilde{\varphi}_{\log \lambda + \log t}$$
noting $g_0=f_1$. using (1.4),(1.5) and definition of the flow we get the result.$\square$

In view of Theorem2.7 and 3.1 we get

{\bf 3.5 Corollary: } {\em A $\ma$ process $\str$ with $X_0=0$ is s.i. $H$-ss iff there exist spectral representation $\{f_t,t\in \re\}$ and two groups of isometries $\{P^{\tau}\}$ and $\{R^{\lambda}\}$ of $L^{\alpha}\bp$ satisfying conditions of Theorem 2.7 and 3.1.}

We do not know, in general, for which $\{P^{\tau}\}$ and $\{R^{\lambda}\}$ can we get a solution $f_t$. We know that these groups are related to flows and cocycles. Suppose the cocycles are constant then we can give flow $\{\varphi_t,t\in\re\}$ and $\{\tilde{\varphi}_{\lambda},\lambda\in\re\}$ so that we can get solutions using corollary 2.9 and Theorem 3.2. However, the following theorem gives most of known examples of $\ma$ s.i,ss processes.

{\bf 3.6 Theorem: } {\em Let $({\cal S},\mu)$ be a standard Borel space and $0<\alpha<2$. Assume

(i) There exists a $\mu$-preserving flow (i.e. $\mu\circ\varphi_t^{-1}=\mu$), $\{\varphi_t,t\in\re\},\;\;P^{\tau}g=g\circ\varphi_{\tau}\\ \tau\in\re$.

(ii) There exists a $g_0\in L^{\alpha}({\cal S},\mu)$ such that $f_t=g_0\circ\varphi_t-g_0\in L^{\alpha}({\cal S},\mu)$ for all $t\in\re$.

(iii) There exists a non-singular flow $\{\tilde{\varphi}_{\lambda},\lambda >0\}$ such that for $r_1\in \re$ constant
$$\frac{d\mu\circ\tilde{\varphi}_t^{-1}}{d\mu}=\lambda^{r_1\alpha}\qquad\mu-a.e.\qquad \forall\lambda>0$$

(iv) $f_{\lambda t}=\lambda^{r_2}f_t\circ\tilde{\varphi}_{\lambda}$ for $t>0,\lambda>0$.\\
Then $X_t=\int_sf_t(s)dM(s)$ is $\ma$ $H$-sss.i. with $H=r_1+r_2$}

The above $P^{\tau}$ and $R^{\lambda}$ satisfy conditions of corollary 3.5 giving the result.

In the remaining parts, we discuss special cases of dissipative and conservative $\ma$ s.i. $H$-ss processes. As stated in Theorem 2.11, an example of a dissipative $\ma$ s.i. process $\str$ with $X_0=0$ is of the form
$$X_t=\int_{W\times \re}[f(x,t-s)-f(x,-s)]N(dx,ds)$$
where $N$ is $\ma$ independently scattered measure on $W\times \re$ with control measure $Q\otimes\lambda$ and $f:W\times \re\fn\re$ is measurable function such that for $t\in\re$
$$\int_{W\times \re}|f(x,t-s)-f(x,s)|^{\alpha}Q(dx)\lambda(ds)<\infty \leqno(3.6')$$
and LHS above tends to 0 as $t\to 0$. Under these conditions, the process $\str$ is well-defined, and has $\ma$ stationary increments, as
\begin{eqnarray*}
X_{t+b}-X_b & = & \int_{{\cal X}\times\re}[f(x,t+b-s)-f(x,b-s)]M(dx,ds)\\
& = & \int_{{\cal X}\times\re}[f(x,t-s)-f(x,-s)]M(dx,d(s-b)),\qquad t\in \re
\end{eqnarray*}
and for any $b\in\re$ the measure $M(dx,ds)$ and $M(dx,d(s-b))$ have the same distribution.

For any $T\in \re$ introduce the `increment process'
$$Y_T(t):=X_{t+T}-X_t = \int_{{\cal X}\times\re}f_T(x,t-s)M(dx,ds)$$
where
$$f_T(x,t)=f(x,T+t)-f(x,t).$$
The process$\{Y_T(t),t\in\re\}$ is a $\ma$ generalized moving average (\sug). Let $\nu_T$ denote its corresponding measure on $L^{\alpha}(\re)$ defined by
$$\nu_T(A)=\int_{L^{\alpha}(\re)}1_A\left(\frac{y(\cdot)}{\|y\|_{\alpha}}\right)\|y\|_{\alpha}^{\alpha}(Q\circ f_T^{-1})(dy).$$
For any $c>0$, let $(S_cy)(t)=y(ct),\;\;t\in\re$, be the scaling transformation, which is a one-to-one mapping of $L^{\alpha}(\re)$ onto itself.

{\bf 3.7 Theorem: } {\em The $\ma$ s.i. process $X=\str$ above is $H$-ss if and only if for all $c>0$ and $T\in\re$,
$$c^{1-\alpha H}\nu_{cT}\circ S_{c^{-1}}\circ \pi^{-1}=\nu_T\circ\pi^{-1}.\leqno(3.8)$$}

{\it Proof: } Note that $X$ is $H$-ss if and only if  for all $c>0$ and $T\in\re$, with $Y(t)=Y_t$
$$\{Y_{cT}(ct):t\in\re\}\sm\{c^HY_T(t):t\in\re\}.\leqno(3.9)$$
Indeed (3.9) follows immediately from the self-similarity of $X$. Conversely, let us show that (3.9) implies that $X$ is $H$-ss. As $X$ is stochastically continuous and $X_0=0$, it suffices to verify that
$$\sum_t b(t)X_{ct}\sm c^H\sum_t b(t)X_t \leqno(3.10)$$
holds for all $c>0$ and any $b(\cdot)$ such that $b(0)=0$ for all but finitely many $t\in L_n$ where $L_n=\{(k_12^{-n},\cdots,k_d2^{-n});\;k_1,\cdots,k_d\in \z\}$ with $\sum_tb(t)=0\;(n=1,2,\cdots)$. But for such $b(\cdot)$, (3.10) can be rewritten as, with $X(t)=X_t$,
$$\sum_t a(t)\left[X(ct+c2^{-n})-X(ct)\right]\sm c^H \sum_t a(t))\left[X(t+2^{-n})-X(t)\right],\leqno(3.11)$$
where $a(t)=\sum_{s\leq t,s\in L_n}b(s)$ if $t\in L_n$, and =0 if otherwise. Clearly (3.11) follows from (3.9) with $T=2^{-n}$.

Introduce $\tilde{Y}_c(t)=c^{1/\alpha}\int_{W\times\re}S_cf_{cT}(x,t-s)M(dx,ds), \;t\in\re$. Then
$$\{Y_{cT}(ct):t\in\re\}\sm\{\tilde{Y}_c(t):t\in\re\}$$
and from Theorem 1 of \sug we have that (3.9) is equivalent to
$$\tilde{\nu}_c\circ\pi^{-1}=c^{\alpha H}\nu_T\circ\pi^{-1}, \leqno(3.12)$$
where $\nu_T$ is given by (3.7) and
\begin{eqnarray*}
\tilde{\nu}(A) & = & c \int_W1_A\left(\frac{S_cf_{cT}(x,\cdot)}{\|S_cf_{cT}(x,\cdot)\|_{\alpha}}\right)\|S_cf_{cT}(x,\cdot)\|_{\alpha}^{\alpha}Q(dx)\\
& = & c\nu_{cT}(S_c^{-1}A),
\end{eqnarray*}
i.e. $\tilde{\nu}=c\nu_{cT}\circ S_c^{-1}=c\nu_{cT}\circ S_{c^{-1}}$. Substituting the above into (3.12) we have (3.8), or the statement of the theorem. $\square$

{\bf 3.13 Corollary } {\em Assume that for every $c>0$ there is a one-to-one `onto' transformation $\rho_c:W\fn W$ such that for every $t\in \re$,
$$f_{ct}(\rho_cs,cs)=c^{\beta_1}f_1(x,s)\qquad Q\otimes Leb-a.e.\leqno(3.14)$$
and
$$Q\circ\rho_c=c^{\beta_2}Q,\leqno(3.15)$$
where $\beta_1,\beta_2\in\re$ are independent of $c>0$. Then $\str$ is $H$-ss with $H=(\alpha\beta_1+\beta_2+1)/\alpha$.}

{\it Proof: } By the definition of $\nu_t$ and (3.14),(3.15),
\begin{eqnarray*}
\nu_{cT}\circ S_{c^{-1}}(A) & = & \int_W1_A\left(\frac{f_{ct}(x,c\cdot)}{\|f_{ct}(x,c\cdot)\|_{\alpha}}\right)\|f_{ct}(x,c\cdot)\|_{\alpha}^{\alpha}Q(dx)\\
& = & \int_W1_A\left(\frac{f_{ct}(\rho_cx,c\cdot)}{\|f_{ct}(\rho_cx,c\cdot)\|_{\alpha}}\right)\|f_{ct}(\rho_cx,c\cdot)\|_{\alpha}^{\alpha}Q(dx)\\
& = & c^{\alpha\beta_1+\beta_2}\nu_t(A).
\end{eqnarray*}
Now apply Theorem 3.7. $\square$

Bellow we consider some concrete cases of the above Corollary. As we shall see, the corresponding $H$-ss processes are numerous and considerably extend the known cases of stable $H$-sssi processes.

{\bf Linear fractional stable motion. } Let $W=\{1\},\;Q=1$ i.e.
$$X_t=\int_{\re}[f(t-s)-f(-s)]M(ds),\qquad t\in\re,\leqno(3.16)$$
where for any $t\in \re$, $f(t-\cdot)-f(-\cdot)\in L^{\alpha}(\re)$, Equation (3.8) in this case is equivalent (see Corollary2,Section 2 of \sug) to
$$c^{(1-\alpha H)/\alpha}f_{cT}(c\cdot)\sim f_T(\cdot),$$
or
$$c^{(1-\alpha H)/\alpha}[f(c(T+s))-f(cs)]=\epsilon [f(T+s+u)-f(s+u)]\qquad ds-a.e.\leqno(3.17)$$
for any $T\in\re,\;c>0$ and some $\epsilon=\epsilon(c,T)=\pm1$ and $u=u(c,T)\in\re$. In particular, if $\epsilon=1$ and $u=0$, then (3.17) coincides with (3.14) (where $t=T$ and $\beta_1=\alpha H-1$), or
$$f(c(T+s))-f(cs)=c^{H-1/\alpha}[f(T+s)-f(s)]\qquad ds-a.e.\leqno(3.18)$$
for any $T\in\re$ and $c>0$. If, in addition, $\hn$ and (3.18) holds for all $s\in\re$, then for some $b_1,b_2\in\re$,
$$f(t)=b_1t_+^{H-1/\alpha}+b_2t_-^{H-1/\alpha}\leqno(3.19)$$
according to Vervaat (1987), and (3.16) becomes the {\em linear fractional stable motion}
$$\begin{array}{lll}
X_t & = & \int_{\re}\{b_1[(t-s)_+^{H-1/\alpha}-(-s)_+^{H-1/\alpha}]\\
&&\qquad +b_2[(t-s)_-^{H-1/\alpha}-(-s)_-^{H-1/\alpha}]\}M(ds), \qquad t\in\re,
\end{array}\leqno(3.20)$$
which is well-defined for $0<H<1,H\neq 1/\alpha$, and is $H$-sssi. The cases $b_1=0$ and $b_0=0$ correspond to the right and the left linear fractional stable motion, respectively. It seems unlikely that (3.18) or even (3.19) in the case $H\neq 1/\alpha$ could have other solutions satisfying $f(t-\cdot)-f(-\cdot)\in L^{\alpha}(\re)$ except (3.19). According to (CMS(1992),Theorem 3), if $f(\cdot)$ is locally integrable, and $H\neq 1/\alpha$, then any non-degenerate $H$-sssi process of the form (3.16) is linear fractional stable motion, the local integrability condition being probably superfluous (see the remark at the end of section 2 in CMS (1992)).

Finally, in the case $H= 1/\alpha$, the solution of Eq (3.18) are
$$f(t)=b_11_{\re_+}(t)+b_21_{\re_-}(t)\qquad (b_1,b_2 \in\re)$$
and
$$f(t)=c\log|t| \qquad (c\in\re),$$
yielding the {\em linear $\ma$ motion}
$$X_t=b_1M(t)+b_2M(-t),\qquad t\in \re,$$
 $M(t):=\int_0^tM(ds),\;t\in\re$, and the {\em log-fractional $\ma$ motion}
$$X_t=c\int_{\re}\log|\frac{t-s}{s}|M(ds),\qquad t\in\re,$$
respectively, which are $H$-sssi with $H= 1/\alpha$; see Vervaart (1987).

{\bf Mixed linear fractional stable motion. } Let ${\cal S}=\re^2=\{b=(b_1,b_2):b_1,b_2\in\re\}$, and let $f(t)=f(b,t)$ be given (3.19), where $\hn$. Condition ($3.6'$) is now equivalent to
$$0<H<1 \leqno(3.21)$$
and
$$\int_{\re^2}|b|_{\alpha}^{\alpha}Q(db)<\infty\leqno(3.22),$$
where $|b|_{\alpha}=(|b_1|^{\alpha}+|b_2|^{\alpha})^{1/\alpha}$. Indeed (3.21) and (3.22) obviously imply ($3.6'$), while on the other hand, if $t>0$, then
\begin{eqnarray*}
\lefteqn{\int_{-\infty}^{\infty}|f(b,t-s)-f(b,-s)|^{\alpha}ds}\\
& & \geq  |b_1|^{\alpha}\int_{-\infty}^0|(t-s)_+^{\ha}-(-s)_+^{\ha}|^{\alpha}ds +|b_2|^{\alpha}\int_t^{\infty}|(t-s)_-^{\ha}-(-s)_-^{\ha}|^{\alpha}ds\\
& & = |b|_{\alpha}^{\alpha}\int_0^{\infty}|(t+s)^{\ha}-s^{\ha}|^{\alpha}ds,
\end{eqnarray*}
with the last integral convergent if and only if (3.21) is satisfied. Hence ($3.6'$) implies (3.21) and (3.22). Let $\rho_cb=b$ for any $c>0$, $b\in\re$; then the conditions of the above corollary are satisfied with $\beta_1=\ha$ and $\beta_2=0$. Consequently the $\ma$ process
$$\begin{array}{lll}
X_t & = & \int_{\re^2\times \re}\{b_1[(t-s)_+^{\ha}-(-s)_+^{\ha}]\\
&& \qquad+b_2[(t-s)_-^{\ha}-(-s)_-^{\ha}]\}M(db,ds), \qquad t\in\re,
\end{array}\leqno(3.23)$$
is well-defined for any $0<\alpha<2,\;0<H<1,\;\hn$ and any $\sigma$-finite measure $Q$ on $\re^2$ satisfying (3.22), and is $H$-sssi according to Corollary 3.13. We shall call (3.23) a {\em mixed linear fractional stable motion with mixing measure} $Q$. Of course, the distribution of (3.23) is not a mixture of the distribution of linear fractional stable motions.

A natural question is how large is the class of ``mixtures'' (3.23), and when the process (3.23) corresponding to different mixing measure are distinct. The answer to the last question is given in Proposition 3.24 bellow. As it turns out, the ``mixture'' form a rather large class, with the linear fractional stable motions corresponding to very special $Q$'s: see Corollary 3.36.

For any measure $\mu$ on a linear topological space $\cal Y$, denoted by $\mu^{(sym})$ the symmetric measure $\mu^{(sym})(dy)=\mu(dy)+\mu(-dy)$.

{\bf 3.24 Proposition } {\em There is a one-to-one correspondence between the distribution of a mixed linear fractional stable motion with mixing measure $Q$, and the measure $\mu_Q^{(sym)}$ on $B_{\alpha}=\{b\in \re^2:|b|_{\alpha}=1\}$, where
$$\mu_Q(do)=\int_0^{\infty}r^{\alpha}Q(dr\times do),\qquad o\in B_{\alpha}.\leqno(3.25)$$}

{\it Proof: } It is clear from the characteristic function of the process (3.23) that $\mu_Q^{(sym)}$ completely determines its distribution.

Below we prove the converse, namely, that if $\{X_t^i: t\in\re\}$ are two mixed linear fractional stable motions with the corresponding mixing measure $Q_i, \;i=1,2$, such that
$$\{X_t^1: t\in\re\}\sm \{X_t^2: t\in\re\},\leqno(3.26)$$
then
$$\mu_{Q_1}^{(sym)}=\mu_{Q_2}^{(sym)}.\leqno(3.27)$$
Consider the ``increment process''
$$Y_t^i:=X_{t+1}^i-X_t^i=\int_{\re^2 \times \re}F(b,t-s)M_i(db,ds),\qquad t\in\re,i=1,2,\leqno(3.28)$$
where
$$\begin{array}{lll}
F(b,t) & = & f(b,t+1)-f(b,t)\\
& = & b_1[(t+1)_+^{\ha}-t_+^{\ha}]+b_2[(t+1)_-^{\ha}-t_-^{\ha}],
\end{array}\leqno(3.29)$$
 $b=(b_1,b_2)\in\re^2$. Then (3.26) implies $\{Y_t^1:t\in\re\}\sm \{Y_t^2:t\in\re\}$ and consequently
$$\nu_1\circ\pi^{-1}=\nu_2\circ\pi^{-1}\leqno(3.30)$$
according to Theorem 1 of \sug, where
$$\nu_i(A)=\int_{B_{\alpha}}1_A\left(\frac{F(o,\cdot)}{\|F(o,\cdot)\|_{\alpha}}\right)
\|F(o,\cdot)\|_{\alpha}^{\alpha}\mu_{Q_i}(do)\leqno(3.31)$$
is a finite measure on the unit sphere ${\cal S}_{\alpha}\subset L^{\alpha}(\re),\;i=1,2$, and the map $G: B_{\alpha}\fn L^{\alpha}(\re)$ by
$$G(o,t):=\frac{F(o,\cdot)}{\|F(o,\cdot)\|_{\alpha}}.\leqno(3.32)$$
Using (3.31), we can rewrite (3.30) as $\mu_1\circ G^{-1}\circ\pi^{-1}=\mu_2\circ G^{-1}\circ\pi^{-1}$, or
$$\mu_1\circ (G\circ\pi)^{-1}=\mu_2\circ (G\circ\pi)^{-1}.\leqno(3.33)$$
Let us prove that (3.33) implies
$$\mu_1^{(syn)}=\mu_2^{(sym)},\leqno(3.34)$$
from which (3.27) easily follows. write $B_{\alpha}=B_{\alpha}^+\cup B_{\alpha}^-\;B_{\alpha}^+=-B_{\alpha}^-,\;;B_{\alpha}^+=\{(o_1,o_2)\in B_{\alpha}: \arctan(o_2/o_1)\in [0,\pi]\}$. Then (3.34) follows from (3.33) if we show that the map $\pi\circ G: B_{\alpha}^+\fn L^{\alpha}(\re)/\!\!\sim$ is continuous and one-to-one; indeed, in this case for any compact set $A\subset B_{\alpha}^+$ the set $B=(\pi\circ G)(A)$ is compact in $L^{\alpha}(\re)/\!\!\sim$ hence also Borel, and such that $(\pi\circ G)^{-1}B=A\cup (-A)$. Consequently from (3.32) we obtain
$$\mu_1^{(sym)}(A)=\mu_1\circ(\pi\circ G)^{-1}B=\mu_2\circ(\pi\circ G)^{-1}B=\mu_2^{(sym)}(A)$$
for any compact set  $A\subset B_{\alpha}^+$ which implies (3.34).

To end the proof of theorem, it remains to verify the continuity and invertibality of the map $\pi\circ G: B_{\alpha}^+\fn L^{\alpha}(\re)/\!\!\sim$. The first property follows from the continuity of\\$G: B_{\alpha}^+\fn L^{\alpha}(\re)$, as $\pi:L^{\alpha}(\re)\fn L^{\alpha}(\re)/\!\!\sim$ is continuous. Note that $F:\re^2\fn L^{\alpha}(\re)$ of (3.33) is continuos, which follows from its linearity or can be verified directly. Moreover, $\|F(o,\cdot)\|_{\alpha}$ nowhere vanishes and is continuous on $B_{\alpha}$, which implies that $\|F(o,\cdot)\|_{\alpha}^{-1}$ is continuous and finally $G(o,\cdot)$ is continuous on $B_{\alpha}$.

Let us first show that $G: B_{\alpha}\fn L^{\alpha}(\re)$ is invertible. Indeed,  $F:\re^2\fn L^{\alpha}(\re)$ is invertible, as $b_1,b_2$ determine the asymptotic of $F(b,t)$ at $+\infty$ and $-\infty$, respectively, i.e. $b_1=(\ha)^{-1}\lim_{t\to +\infty}F(b,t)t^{1+1/\alpha-H}, \;b_2=(\ha)^{-1}\lim_{t\to -\infty}F(b,t)t^{1+1/\alpha-H}$. Therefore $G(o_1,\cdot)=G(o_2,\cdot)\;(o_1,o_2\in B_{\alpha})$ implies $o_1\|F(o_1,\cdot)\|_{\alpha}^{-1}=o_2\|F(o_2,\cdot)\|_{\alpha}^{-1}$ and $o_1=o_2$.

Now let $o_1,o_2 \in B_{\alpha}^+$ be such that $\pi\circ G(o_1)=\pi\circ G(o_2)$ and $o_1\neq o_2$. Hence $\pi(G(o_1,\cdot))=\pi(G(o_2,\cdot))$ which implies $G(o_1,\cdot)\sim G(o_2,\cdot)$ or
$$G(o_1,\cdot)\sim \epsilon G(o_2,\cdot)\leqno(3.35)$$
for some $\epsilon=\pm 1$ and $t\in\re$. Now, if $t=0$ then (3.35) implies $o_1=\epsilon o_2$ by the invertibility and the definition of $G$ which contradicts $o_1\neq o_2,\;o_1,o_2 \in B_{\alpha}^+$. If $t\neq 0$, then (3.35) is false again, e.g. because $F(b,\cdot)$ and $G(b,\cdot)$ are analytic on the intervals $(-\infty,-1),(-1,0),(0,\infty)$ and not analytic at $t=-1$ and $t=0$ unless $b=0$. Hence $o_1=o_2$ which proves the invertibility of $\pi\circ G|_{B_{\alpha}^+}$ and the the proposition. $\square$

The measure $\mu_Q^{(sym)}$ of (3.25) is concentrated at two (symmetric) points of $B_{\alpha}$ if and only if $Q$ is concentrated on a ray (i.e. a line in $\re^2$ going through origin). Together with Proposition 3.24 this implies

{\bf 3.36 Corollary } {\em The distribution of a mixed linear fractional stable motion coincides with the distribution of a linear fractional stable motion if and only if the mixing measure $Q$ is concentrated on a ray.}

If we do not distinguish between processes up to constant multipliers, then two mixed linear fractional stable motions will be distinct if and only if their $\mu_Q^{(sym)}$ measures are not multiples of each other; and as a Corollary we obtain the characterization of distinct linear fractional stable motions by means of distinct lines through the origin of the parameter plane $(b_1,b_2)$ given in Theorem 3.1 of Cambanis and Maijima (1989).

{\bf Mixed truncated (left) fractional stable motion. } Consider the following class of $\ma$ processes which for all $a,b\in\re$ are given by
$$X_t^{a,b}=\int_{\re_+\times\re}[(t-s)_+^a\wedge p^a-(-s)_+^a\wedge p^a]M(dp,ds),\qquad t\in\re\leqno(3.37)$$
where $M$ is a $\ma$ independently scattered random measure on $\re_+\times\re$ with control measure $Q\otimes Leb$, and
$$Q(dp)=p^{-1-b}dp,\qquad p\in\re_+.\leqno(3.38)$$
We call (3.37) a {\em mixed truncated left fractional stable motion with parameters $a,b$.}

{\bf 3.39 Proposition } {\em The process $\{X_t^{a,b}:t\in\re\}$ is well-defined if and only if either
$$0\vee (\alpha a-\alpha+1)<b<\alpha a \leqno(3.40)$$
or
$$\alpha a<b<0 \wedge (\alpha a+1) \leqno(3.41)$$
holds, and in this case it is $H$-sssi with
$$H=\frac{\alpha a-b+1}{\alpha},\leqno(3.42)$$
and the distribution of $\{X_t^{a,b}:t\in\re\}$ corresponding to different pairs $(a,b)$ are distinct.}

{\it Proof: } Consider the integral
$$I(a,b)=\int_0^{\infty}\int_{-\infty}^{\infty}|(t-s)_+^a\wedge p^a-(-s)_+^a\wedge p^a|^{\alpha}p^{-1-b}dpds.$$
Let $a<0$. Then for $t>0$,
\begin{eqnarray*}
I(a,b) & = & \int_0^{\infty}ds\int_0^sdp \;p^{-1-b}|(t+s)^a-s^a|^{\alpha}\\
&& \qquad +\int_0^{\infty}ds\int_s^{t+s}dp \;p^{-1-b}|(t+s)^a-s^a|^{\alpha}\\
&& \qquad +\int_0^tds\int_0^s dp\;p^{-1-b}s^{\alpha a}+\int_0^tds\int_0^{\infty}dp\;p^{-1-b}p^{\alpha a}=:\sum_{i=1}^4I_i.
\end{eqnarray*}
Here, $I_3$ and $I_4$ are finite if and only if $b<0,\alpha a+1>b$ and $b>\alpha a$, i.e. condition (3.41) is nacessary. It is easy to check that $I_1$ and $I_2$ are also finite under (3.41). In the case $a>0$ we can write similarly $I(a,b)=\sum_{i=1}^4I_i$, where
\begin{eqnarray*}
I_1 & = & \int_0^{\infty}ds\int_s^{t+s}dp \;p^{-1-b}|s^a-p^a|^{\alpha},\\
I_2 & = & \int_0^{\infty}ds\int_{t+s}^{\infty}dp \;p^{-1-b}|(t+s)^a-s^a|^{\alpha},\\
I_3+I_4 & = & \int_0^tds\int_0^s dp\;p^{-1-b}p^{\alpha a}+\int_0^tds\int_0^{\infty}dp\;p^{-1-b}s^{\alpha a}.
\end{eqnarray*}
Here $I_3+I_4<\infty$ if and only if $0<b<\alpha a$. But then $I_2<\infty$ iff $b>1+\alpha a-\alpha$ and $I_1<\infty$ if $b>\alpha a-\alpha$, i.e. (3.40) is necessary and sufficient. Note also that $I(a,b)=\infty$ if $a>0$ and $0<\alpha\leq 1$. Finally, for $a=0$ the integral $I(a,b)$ clearly diverges.

Let $\rho_cp=cp,\;p\in \re_+$. Then conditions of Corollary 3.13 are satisfied for $\beta_1=a$ and $\beta_2=-b$, and therefore $\{X_t^{a,b}:t\in\re\}$ is $H$-sssi with $H$ given by (3.42).

To prove the last statement of the proposition, consider the ``increment process''
$$Y_{a,b}(t)= X_{t+1}^{a,b}-X_t^{a,b}=\int_{\re_+^2\times\re}F(p,t-s)M(dp,ds),\qquad t\in\re ,$$
where $F(p,t)=(t+1)_+^a\wedge p^a-t_+^a\wedge p^a$. The corresponding measure $\nu_{a,b}$ is concentrated on the set $K_a=\{cF(p,\cdot):c\in\re,p\in\re_+\}$. As $K_a\cap U_tK_a=\emptyset$ for $t\neq 0$ and $K_{a'}\cap K_{a''}=\emptyset$ for $a'\neq a''$ we see that $\nu_{a',b'}\circ\pi^{-1}\neq \nu_{a'',b''}\circ\pi^{-1}$ for any $(a',b'),(a'',b'')$ such that $a'\neq a''$. Therefore $Y_{a',b'}\sm Y_{a'',b''}$ and $X^{a',b'}\stackrel{d}{\neq} X^{a'',b''}$ if $a'\neq a''$. In the case $a'=a''\equiv a$ and $b' \neq b''$ we have again $X^{a',b'}\stackrel{d}{\neq} X^{a'',b''}$ as $X^{a',b'}$ and $X^{a'',b''}$ are $H$-ss with $H=H'=(\alpha -b'+1)/\alpha$ and $H=H''=(\alpha -b''+1)/\alpha$,  respectively, and $H'\neq H''$. Proposition (3.39) is proved.$\square$

{\bf Remark } One can show, moreover, that the process $X_{a,b}$ of (3.37) are different in the sense of distribution from any mixed linear fractional stable motion. It is also possible to consider more general $H$-sssi processes of this type, e.g.
\begin{eqnarray*}
X_t^{a,b,c} & = & \int_{\re_+^2\times\re}\left\{c_1[(t-s)_+^a\wedge p_1^a-(-s)_+^a\wedge p_1^a]\right.\\
&& \qquad +\left. c_2[(t-s)_-^a\wedge p_2^a-(-s)_-^a\wedge p_2^a]\right\}M(dp,ds),\qquad i\in\re,
\end{eqnarray*}
$c=(c_1,c_2)\in\re^2$, with $Q(dp)=p_1^{-1-b}p_2^{-1-b}dp_1dp_2,\;p=(p_1,p_2)\in\re_+^2$, and even more general ``truncated mixtures''.

{\bf Chentsov type stable $H$-sssi processes. } Recently Tekenaka (1991) introduced a new class of stable self-similar processes called {\em generalized Chentsov type} which are defined by
$$X_t=N(S(t)),\qquad t\in\re\leqno(3.43)$$
where $N$ is a $\ma$ independently scattered random measure on $\re_+\times \re$ with control measure $Q\otimes \lambda$, $Q(dx)=x^{2-\beta}dx,\;x\in\re_+$,
$$S(t)=C(0)\Delta C(t),$$
 $\Delta$ is symmetric difference of sets, and
$$C(t)=\{(x,s)\in\re_+\times\re :|t-s|<x\},\qquad t\in\re.$$
The process (3.43) is well-defined if and only if $0<\beta<1$ and is $H$-sssi with $H=\beta/\alpha$ (Theorem 4 of Takenaka (1991)). Sato (1991) proved that the distribution of $\str$ is completly determined by its bivariate distribution (Proposition 1 of \sug). In Talenaka and Sato, multiparameter versions of (3.43), or $\ma$ self-similar random fields in $\re^d\;(d\geq 1)$, were studied also.

The process (3.43) is a particular case of Corollary 3.36, and gives a more genuine example of a process of the form we began with where mixing is essential. Indeed (3.43) can be written as
$$X_t=\int_{\re_+\times\re}[1(|t-s|<x)-1(|-s|<x)]M(dx,ds),\qquad t\in\re,\leqno(3.44)$$
where
$$M(dx,ds):=\left\{\begin{array}{ll}
N(dx,ds), & \mbox{if } (x,s)\notin C(0),\\
-N(dx,ds), & \mbox{if } (x,s)\in C(0),
\end{array}\right. $$
and $M$ is again $\ma$ measure on $\re_+\times\re$ having the same distribution as $N$. It is easy to check that for $X_t$ of (3.44), the conditions of corollary 3.36 are satisfied with $\rho_c(x)=cx,\;x\in\re_+,\;\beta_1=0$ and $\beta_2=\beta-1$, hence $\str$ is $H$-sssi with $H=\beta/\alpha \in (0,1/\alpha)$. Note that unlike the remaning examples discussed in this section, the processes (3.43) or (3.44) can be $H$-ss with $H>1$ (provided $0<\alpha <1$).

\section {``Conservative'' $\ma$ siss processes}
Below, we use a similar construction to Section 3 to define a class of conservative $\ma$ ssssi processes. Consider the measure preserving flow
$$\phi_t(s,x)=(s+tx_{mod\,2\pi},x),\qquad t\in \re\leqno(4.1)$$
on $U=(0,2\pi)\times \re_+$ with $\mu=$Leb$\otimes Q$, where $Q$ is any $\sigma$-finite measure on $\re_+$. The flow (4.1) defines the group
$$P_tf=f\circ\phi_t,\qquad t\in \re,\;f\in L^{\alpha}(U)\leqno(4.2)$$
of isomeries of $L^{\alpha}(U)$. The corresponding stationary $\ma$ process
$$X_t=P_tf,\qquad t\in \re\leqno(4.3)$$
will be called a mixed rotating average with the corresponding parameters $(f,Q)$. The above terminology is consistent with \sug and Section 3 and reflects the fact that the transformation $s\fn s+tx_{mod\,2\pi}$ is rotation of the circle $(0,2\pi)$ with
speed $x$. It is easy to show that the flow (4.1) is conservative. The process (4.3) admits the spectral representation
$$X_t=\int_{(0,2\pi)\times\re^2}f(s+tx_{mod\,2\pi},x)M(ds,dx)\leqno(4.4)$$
where $\mu$ is a $\ma$ measure on $(0,2\pi)\times\re^2$ with control measure Leb$\otimes Q$. The case where $Q$ is concentrated at one point $x_0\in\re_+$ corresponds to the periodic process
$$X_t=\int_0^{2\pi}f(s+tx_0)M(ds),\qquad t\in \re\leqno(4.5)$$

{\bf 4.6 Remark } In the case $\alpha=2$, any stationary ergodic Gaussian process admits the representation (4.1). Indeed, let $f(s,x)=\cos(s)$ and $Q(\re_+)<\infty$. Then
\begin{eqnarray*}
r(t)=EX_0X_t & = & \int_{\re_+}dQ(x)\int_0^{2\pi}\cos(s+tx)\cos(s)ds\\
& =& \int_{\re_+}\cos(tx)dF(x),
\end{eqnarray*}
where $dF(x)=\mbox{const.}\,dQ(x)$ is an arbitrary finite measure.

In the sequel, we shall need a criterion to decide when two mixed rotating average $X^1,X^2$ with parameters $(f_1,Q_1),(f_2,Q_2)$ have the same distribution. If $Q_1=Q_2$, then clearly $f_1(s,x)$ and $f_2(s,x)=\epsilon f_1(s+\tau(x),x)$, where $\epsilon=\pm 1$ and $\tau(x)$ is any measurable function $\re_+\fn\re_+$, yield equivalent spectral representation.

Below, we identify any function $f$ on $[0,2\pi)\times \re_+$ with its periodic extension to $\re\times\re_+$, i.e. $f(s,x)=f(s+2k\pi,x),\; \forall (s,x)\in[0,2\pi)\times \re_+,\;\forall k\in\z$.

{\bf 4.7 Definition } a function $f$ on $[0,2\pi)\times \re_+$ is said periodically minimal if $2\pi$ is the minimal period of $f(\cdot,x)$ for every $x\in\re_+$.

{\bf 4.8 Proposition } {\em Let $f=f(s,x)$ be a measurable function on $[0,2\pi)\times \re_+$. The representation $\{f(s+tx,x),t\in\re\}$ is minimal in the sense of Hardin (1982) iff there exists a periodically minimal function $g=g(s,x)$ such that $f(s,x)=g(s,x)\;Leb\otimes Q$-a.e.}

 {\it Proof: } To prove the sufficient part, it suffices to show that the map $(s,x)\fn\{g(s+tx,x):t\in\re\}$ is one-to-one. Indeed, let $(s_1,x_1),(s_2,x_2)\in [0,2\pi)\times \re_+,\;s_1<s_2$, be such that
$$g(s_1+tx_1,x_1)=g(s_2+tx_2,x_2)\qquad \forall t\in\re.$$
As both sides are periodic in $t$, so their periods have to be equal, which yields $x_1=x_2\equiv x$. Thus,
$$g(s_1+tx_,x_)=g(s_2+tx_,x)\qquad \forall t\in\re.$$
or
$$g(s,x)=g(s+s_2-s_1,x)\qquad\forall s\in\re.$$
i.e., $g(\cdot,x)$ is periodic with the period $0<s_2-s_1<2\pi$, which contradicts the periodic minimality assumption of $g$.

{\bf 4.9 Theorem }{\em Let $X^1,X^2$ be two mixed rotating averages with parameters $(f_1,Q_1)$, $(f_2,Q_2)$ respectively, where $f_1,f_2$ are periodically minimal. Then $X^1\sm X^2$ iff

(i) $Q_1\sim Q_2$,

(ii) there exists measurable function $\tau:\re_+\fn\re$ and $h$ so that,
$$Uf(s,x)=h(x)f(s+\tau(x),x)$$
extends to a linear isometry $L^{\alpha}(Leb\times Q)\fn L^{\alpha}(Leb\times Q)$, in particular,
$$\int_{\re_+}|h(x)|^{\alpha}|f(x)|^{\alpha}dQ_2=\int_{\re_+}|f(x)|^{\alpha}dQ_1\leqno(4.10)$$
for any $f\in L^{\alpha}(Q_1)$.}

Next, we consider ``conservative'' $\ma$ si processes of the form
$$X_t=\int_{[0,2\pi)\times \re_+} f(t,s,x)M(ds,dx)\leqno(4.11)$$
where
$$f(t,s,x)=g(s+tx)-g(s)=g\circ\phi_t(s,x)-g(s)\leqno(4.12)$$
and $g(\cdot)$ is a measurable function on $[0,2\pi)$ which is periodically extended to $\re$. We assume $f(t,\cdot)\in L^{\alpha}(Leb\times Q), \;\forall t\in\re$, and $f(t,.)\to 0$ in $L^{\alpha}(Leb\times Q)$ as $t\to 0$, i.e., that $X_t$ of (4.11) is stochastically continuous. This last assumption implies in particular that $g\in L^{\alpha}(Leb)$. Indeed, as $\|f(t,\cdot)\|_{L^{\alpha}(Leb\times Q)}=:F(t)$ is continuous so
\begin{eqnarray*}
\infty & > & \int_0^{2\pi}F(t)dt\geq \int_0^{2\pi}\int_{[1,\infty)}Q(dx)\int_0^{2\pi}dt|g(s+tx)-g(s)|^{\alpha}\\
& \geq & \int_0^{2\pi}ds\int_{[1,\infty)}x^{-1}Q(dx) \int_0^{2\pi}dt|g(s+t)-g(s)|^{\alpha}.
\end{eqnarray*}
Hence there is $s_0\in (0,2\pi)$ such that $ \int_0^{2\pi}|g(s_0+t)-g(s_0)|^{\alpha}dt=\int_0^{2\pi}|g(t)-g(s_0)|^{\alpha}dt<\infty$, or $g\in L^{\alpha}(Leb)$.

{\bf Proposition 4.13 } {\em Assume
$$Q(dx)=x^{-1-\beta}\qquad x\in\re_+,\leqno(4.14)$$
 $g\in L^{\alpha}(Leb)$ and $\int_0^{2\pi}|g(s+t)-g(s)|^{\alpha}ds=O(t^{-r})\;(t\to 0)$ for some $r>\beta>0$. Then $X_t$ of (4.11) is well-defined and $H$-ss with $H=\beta/\alpha$.}

Proof reduces to easy verification of conditions of Theorem 3.6 with the flow $\phi_t$ given by (4.1) and $\tilde{\varphi}_{\lambda}(s,x)=(s,\lambda x),\; (s,x)\in [0,2\pi)\times \re_+,\;\lambda >0$.

Finally, we address the question when two ss processes described in the last Proposition coincide in distribution.

{\bf Theorem 4.15 } {\em let $X_t^1, X_t^1$ be two $\ma$ sssi processes of Proposition 4.13 and corresponding to parameters $(g_1,\beta_1),(g_2,\beta_2)$, respectively. Assume that $g_1,g_2 \in L^{\alpha}(Leb)$ are periodically minimal. Then $X^1\sm X^2$ if and only if $\beta_1=\beta_2$ and there exists $\epsilon =\pm\;c\in\re$, and $\tau\in [0,2\pi)$ such that
$$g_2(\cdot)=\epsilon g_1(\cdot +\tau)+c.\leqno(4.16)$$}

{\it Proof: } The sufficiency part being rather obvious, we shall prove only the necessity. Note first that $X^1\sm X^2$ implies $H_1=H_2$, and therefore $\beta_1=\beta_2$ or $Q_1=Q_2$.

Let $X_t^i=\int_{[0,2\pi)\times \re_+} f_i(t,s,x)dM,\; f_i(t,s,x)=g_i(s+tx)-g_i(s),\;i=1,2$, then $X^1\sm X^2$ implies
$$Y_1\sm Y_2,\leqno(4.17)$$
where
$$Y_i(t)=X_{t+1}^i-X_t^i=\int_{[0,2\pi)\times \re_+} f_i(1,\cdot)\circ\phi_tdM,\leqno(4.18)$$
are mixed rotating averages discussed above. Next, in order to be able to aply Theorem 4.9, we need to show that the fact that $g_1(\cdot),g_2(\cdot)$ are periodically minimal implies that $f_1(1,\cdot),f_2(1,\cdot)$ are Leb$\otimes Q$-a.e. equal to periodically minimal functions. Assume that this is not the case; then there exists a set $A\subset\re_+$ with $Q(A)>0$ (or Leb($A)>0$) and an integer $k>1$ such that for any $x\in A,\;g(x+s)-g(s)$ is a.e. equal to a periodic function with period $2\pi/k$. Next, as linear combination of periodic are periodic, we infer that $g(x_2-x_1+s)-g(s)$ is periodic with period $2\pi/k$ for any $x_1,x_2\in A$. But the set $A$ has a density point, so $x_2-x_1$ can be arbitrary small.

Assume for simplicity that $g(\cdot)$ is a continuously differentiable then the above argument implies that $g'(\cdot)$ is periodic with period $2\pi/k$. But this implies that $g(\cdot)$ is periodic with period $2\pi/k$, which is contradiction. Indeed,
\begin{eqnarray*}
0 &=& g(s+2\pi)-g(s)=\int_s^{s+2\pi}g'(u)du\\
&=& \sum_{j=0}^{k-1}\int_{s+\frac{2\pi}{k}j}^{s+\frac{2\pi}{k}(j+1)}g'(u)du\\
&=& k\int_s^{s+\frac{2\pi}{k}}g'(u)du\\
&=& k(g(s+\frac{2\pi}{k})-g(s)),
\end{eqnarray*}
or $g(s+\frac{2\pi}{k})=g(s)$.

Thus, we can apply Theorem 4.9 to $Y_i(t)$ of (4.18), which implies
$$g_2(s+x)-g_2(s)=\epsilon (g_1(s+x+\tau(x))-g_1(s+\tau(x)))\leqno(4.19)$$
for some $\epsilon=\pm 1,\;\tau(x)$ a measurable function, and a.e. $x>0$.
In particular, (4.19) holds for some $x$ arbitary small.

Assume again for simplicity that $g_1(\cdot),g_2(\cdot)$ are continuously differentiable, and that $\tau(x)$ is continuous at $x=0$. Then (4.17) implies
$$g_2'(s)=\epsilon g_1'(s+\tau),\leqno(4.20)$$
Where $\tau=\tau(0)$. Consequently,
\begin{eqnarray*}
g_2(s) & = & g_2(0)+\int_0^sg_2'(u)du\\
&=& g_2(0)+\epsilon \int_0^sg_1'(u+\tau)du\\
&=& \epsilon g_1(s+\tau)+c,
\end{eqnarray*}
where $c=g_2(0)-g_1(\tau)$.

\newpage
\begin{center}
{APPENDIX: ON RIEMANN INTEGRATION OF CURVES}
\end{center}

Consider a function $f:[a,b]\fn E$, where $E$ is a complete normed space. We shall write $f\in{\cal R}_m([a,b],E)$ if for every $\phi \in D^-[a,b]$, Riemann's integral
$$\int_a^b \phi (t) f(t) dt\leqno(A.1)$$
exists. Here $D^-[a,b]$ is the space of all left-continuous and with right-hand limits functions $\phi:[a,b]\fn \re$. $D^-[a,b]$ is a Banach space with the supremum norm in which function of the type
$$\phi = x_01_{\{\alpha\}}+\sum_i x_i1_{(t_{i-1},t_i]}$$
from a dense subspace $(a=t_0< \cdots <t_n=b)$. As a consequence of Banach-Steinhaus theorem we get that for every $f\in{\cal R}_m([a,b],E)$ the map
$$D^-[a,b\ni \phi \fn \int_a^b \phi (t) f(t) dt\in E \leqno(A.2)$$
is continuous. We will need the following form of Fubini's theorem.

{\bf A.3 Lemma } {\em Let $f\in{\cal R}_m([a,b],E)$ and let $F(t)=\int_a^tf(s)ds$. Then $F\in{\cal R}_m([a,b],E)$ and, for every $\phi \in D^-[a,b]$,
$$\int_a^b \phi (t) F(t) dt=\int_a^b\left(\int_t^b \phi (s)ds\right)f(t)dt.$$}

{\it Proof: } Let $\pi:a=t_0< \cdots <t_n=b$ be a partition and $\xi_j\in[t_{j-1},t_j],\;\xi_0=a$. Put $\Delta t_j=t_j-t_{j-1},\;\Delta F_j=F(\xi_j)-F(\xi_{j-1})$, and $|\pi|=\max\Delta t_j$. Then we have
\begin{eqnarray*}
\lefteqn{\sum_{j=1}^n \phi(\xi_j) F(\xi_j) \Delta t_j = \sum_{j=1}^n \phi(\xi_j) \left\{\sum_{i=1}^j \Delta F_i \right\}\Delta t_j} \\
& &  \qquad =\sum_{i=1}^n \left\{\sum_{j=i}^n\phi(\xi_j)\Delta t_j \right\} \Delta F_i = \int_a^b\Phi_{\xi,\pi}(t)f(t)dt
\end{eqnarray*}
where
$$\Phi_{\xi,\pi}=\sum_{i=1}^n \left(\sum_{j=i}^n\phi(\xi_j)\Delta t_j \right)1_{(\xi_{i-1},\xi_i]}.$$
Notice that by the uniform continuity of $\Phi(t):=\int_t^b\phi(s)ds$ and $\Phi(b)=0$ we get $\Phi_{\xi,\pi}\to \Phi$ uniformly as $|\pi|\to 0$. Hence
$$\int_a^b\Phi_{\xi,\pi}(t)f(t)dt\to \int_a^b \Phi (t) f(t) dt$$
as $|\pi|\to 0$, which ends the proof. $\square$

If $E$ is Banach space, then every continuous $f$ belong to ${\cal R}_m([a,b],E)$. However, if $E$ is not locally convex, then there are continuous $E$-valued functions which are not Riemann integrable (see Rolewicz (1984)). In the latter case, even when $f$ is continuous and Riemann integrable, we do not know whether $\phi f$ is Riemann integrable for every\\
$C_{\infty}$-function $\phi$.For this reason we suppose the space ${\cal R}_m([a,b],E)$, which is closed under $D^-$-multipliers and permits some integral calculus. We will not further develo this calculus since, for our purpose, Lemma A.3 is sufficient. We will now consider other conditions guaranteeing the existence of (A.1). We may consider a semivariation
$$A(f,\delta)=\sup\|\sum_j c_j[f(t_j)-f(t_{j-1})]\|$$
where the supremum is taken over $|c_j|\leq\delta ,\;a=t_0< \cdots <t_n=b$, and $n\geq 1$. Then we have the following.

{\bf A.4 Lemma } {\em Suppose thet $\lim_{\delta\to 0}A(f,\delta)=0$. Then $f\in{\cal R}_m([a,b],E)$.}

{\it Proof: } Let $\pi:a=t_0< \cdots <t_n=b$ be a partition and $\xi_j\in[t_{j-1},t_j],\;\xi_0=a$. Put $\Delta t_j=t_j-t_{j-1},\;\Delta f_j=f(\xi_j)-f(\xi_{j-1})$, and $|\pi|=\max\Delta t_j$. For $\phi \in D^-[a,b]$ consider Riemann's sum
\begin{eqnarray*}
\lefteqn{\sum_{\xi,\pi}=\sum_{j=1}^n \phi(\xi_j) f(\xi_j) \Delta t_j = \sum_{j=1}^n \phi(\xi_j) \left\{\sum_{i=1}^j \Delta f_i +f(a)\right\}\Delta t_j} \\
& &  \qquad =\sum_{i=1}^n \left\{\sum_{j=i}^n \phi(\xi_j)\Delta t_j \right\} \Delta f_i + f(a)\sum_{i=1}^n \phi(\xi_j)\Delta t_j.
\end{eqnarray*}
Put $\Phi(x)=\int_x^b\phi(s)ds$ and notice that the continuity of $\Phi$ implies that
$$\sum_{j=i}^n \phi(\xi_j)\Delta t_j=\Phi(\xi_i)+c_i$$
where $|c_i|$ are uniformly small when $|\pi|$ is small. Hence we get
$$\sum_{\xi,\pi}=\sum_{i=1}^n \Phi(\xi_i)\Delta f_i+f(a)\Phi(a)+x_{\xi,\pi}$$
where $\|x_{\xi,\pi}\|\leq A(f,\delta)+\|c_1f(a)\|$ and $\delta=\max |c_i|;\;\|x_{\xi,\pi}\|$ is small when $|\pi|$ is small.

Let now $\pi '$ be another partition of $[a,b]$ and $\xi_j '\in[t_{j-1} ',t_j '],\;k\leq m$. Let $\xi ''$ be the common refinement of $\xi$ and $\xi '$. We have
$$\sum_{\xi,\pi} - \sum_{\xi ',\pi '}=\sum_{i=1}^n[\Phi(\eta_i)-\Phi(\eta_i ')][f(\xi_i '')-f(\xi_{i-1} '')+(|x_{\xi,\pi}-|x_{\xi ',\pi '})$$
where $|\eta_i -\eta_i '|\leq 2|\pi|+2|\pi '|$. Using our assumption on $A(f,\cdot)$ and the uniform continuity of $\Phi$ we get that $\|\sum_{\xi,\pi} - \sum_{\xi ',\pi '}\|$ is arbitrarily small when both $|\pi|$ and $|\pi '|$ are small. This concludes the proof. $\square$

\end{document}